# Parallel Transport Frame in 4–dimensional Euclidean Space $\mathbb{E}^4$.

Fatma GÖKÇELİK, Zehra BOZKURT, İsmail GÖK, F. Nejat EKMEKCİ,
and Yusuf YAYLI

**Abstract :** In this work, we give parallel transport frame of a curve and we introduce the relations between the frame and Frenet frame of the curve in 4–dimensional Euclidean space. The relation which is well known in Euclidean 3–space is generalized for the first time in 4–dimensional Euclidean space. Then we obtain the condition for spherical curves using the parallel transport frame of them. The condition in terms of $\kappa$ and $\tau$ is so complicated but in terms of $k_1$ and $k_2$ is simple. So, parallel transport frame is important to make easy some complicated characterizations. Moreover, we characterize curves whose position vectors lie in their normal, rectifying and osculating planes in 4–dimensional Euclidean space $\mathbb{E}^4$.

## 1. Introduction

The Frenet frame is constructed for the curve of 3-time continuously differentiable non-degenerate curves. But, curvature may vanish at some points on the curve. That is, second derivative of the curve may be zero. In this situation, we need an alternative frame in $\mathbb{E}^3$. Therefore In [1], Bishop defined a new frame for a curve and he called it Bishop frame which is well defined even if the curve has vanishing second derivative in 3–dimensional Euclidean space. The advantages of the Bishop frame and the comparable Bishop frame with the Frenet frame in Euclidean 3–space was given by Bishop [1] and Hanson [9]. In Euclidean 4–space $\mathbb{E}^4$, we have the same problem for a curve like being in Euclidean 3–space. That is, one of the $i-th$ $(1 < i < 4)$ derivative of the curve may be zero. In this situation, we need an alternative frame. So, using the similar idea we consider such curves and construct an alternative frame.

Our work is structured as follows. Firstly, we give parallel transport frame of a curve and we introduce the relations between the frame and Frenet frame of the curve in 4–dimensional Euclidean space. The relation which is well known in Euclidean 3–space is generalized for the first time in 4–dimensional Euclidean space. For construction of parallel transport frame we use the following method.

Let $\alpha(s)$ be a space curve parametrized by arclenght $s$ and its normal vector field be $V(s)$ which is perpendicular to its tangent vector field $T(s)$ said to be relatively parallel vector field if its derivative is tangential along the curve $\alpha(s)$. For a given curve if $T(s)$ is given unique, we can choose any convenient arbitrary basis which consist of relatively parallel vector field $\{M_1(s), M_2(s), M_3(s)\}$ of the frame, they are perpenticular to $T(s)$ at each point.

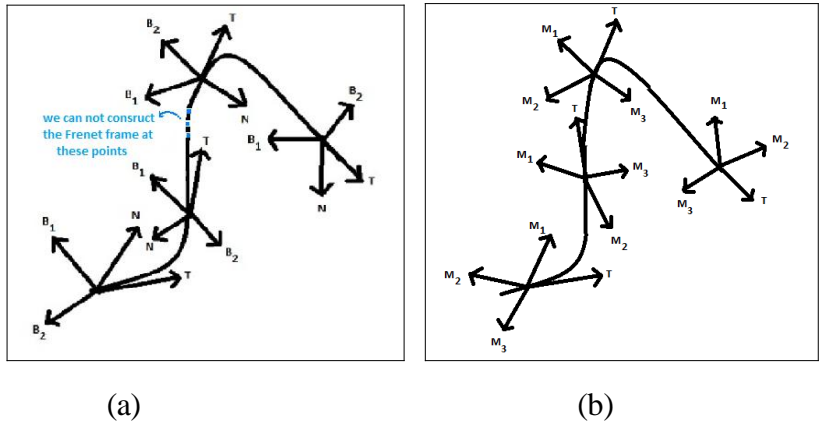

(a)          (b)

Figure 1: Comparing the Frenet and parallel transport frames.

We consider the Frenet frame which is undefined at the some points along the curve and the parallel transport frame is smooth along the curve in Figure 1-(a) and in Figure 1-(b), respectively. Then we obtain the condition for spherical curves using the parallel transport frame of them. The condition in terms of $\kappa$ and $\tau$ is so complicated in Euclidean $3-$space but the condition given by Bishop [1] in terms of $k_1$ and $k_2$ is simple. So, parallel transport frame is important to make easy some complicated characterizations. So, secondly we give a formulae for a curve such as Let $\alpha : I \subset \mathrm{R} \to \mathbb{E}^4$ be a curve with nonzero curvatures $k_i$ $(i=1,2,3)$ according to parallel transport frame in $\mathbb{E}^4$. Then $\alpha$ lies on a sphere if and only if $ak_1 + bk_2 + ck_3 + 1 = 0$ where $a$, $b$, and $c$ are non-zero constants.

Finally, in this work, we give a necessary and sufficient conditions for curves in Euclidean $4-$space $\mathbb{E}^4$ to be normal, rectifying and osculating curves in terms of their parallel transport curvature functions. Normal curves in Minkowski space $\mathbb{E}_1^4$ are defined in [12] as the space curves whose position vector always lies in its normal space $T^\perp$, which represents the orthogonal complement of the tangent vector field of the curve. Rectifying curves are defined in the Euclidean $4-$space as a curve whose position vector always lies in orthogonal complement $M_1^\perp$ of its principal normal vector field $M_1$ [10]. Analogously, osculating curve in Minkowski space-time are defined in [11] as the space curves whose position vector always lies in its osculating space, which represents the orthogonal complement of the second binormal vector field $M_2$ of the parallel transport frame.

## 2. Preliminaries

Let $\alpha : I \subset \mathbb{R} \to \mathbb{E}^4$ be arbitrary curve in the Euclidean $4-$space $\mathbb{E}^4$. Recall that the curve $\alpha$ is parameterized by arclength function $s$ if $\langle \alpha'(s), \alpha'(s) \rangle = 1$ where $\langle , \rangle$ is the inner product of $\mathbb{E}^4$ given by

$$\langle X, Y \rangle = x_1 y_1 + x_2 y_2 + x_3 y_3 + x_4 y_4$$

for each $X = (x_1, x_2, x_3, x_4)$, $Y = (y_1, y_2, y_3, y_4) \in \mathbb{E}^4$. In particular, the norm of a vector $X \in \mathbb{E}^4$ is given by $\|X\| = \sqrt{\langle X, X \rangle}$. Let $\{T, M_1, M_2, M_3\}$ be the moving Frenet frame along the unit speed curve $\alpha$.

Then $T, N, B_1$ and $B_2$ are the tangent, the principal normal, first and second binormal vectors of the curve $\alpha$, respectively. If $\alpha$ is a space curve, then this set of orthogonal unit vectors, known as the Frenet-Serret frame, has the following properties

$$\begin{bmatrix} T' \\ N' \\ B_1' \\ B_2' \end{bmatrix} = \begin{bmatrix} 0 & \kappa & 0 & 0 \\ -\kappa & 0 & \tau & 0 \\ 0 & -\tau & 0 & \sigma \\ 0 & 0 & -\sigma & 0 \end{bmatrix} \begin{bmatrix} T \\ N \\ B_1 \\ B_2 \end{bmatrix}$$

where $\kappa$, $\tau$ and $\sigma$ denote principal curvature functions according to Serret-Frenet frame of the curve $\alpha$, respectively.

We use the tangent vector $T(s)$ and three relatively parallel vector fields $M_1(s)$, $M_2(s)$ and $M_3(s)$ to construct an alternative frame. We call this frame a parallel transport frame along the curve $\alpha$. The reason for the name parallel transport frame is because the normal component of the derivatives of the normal vector field is zero.

We shall call the set $\{T, M_1, M_2, M_3\}$ as parallel transport frame and

$$k_1 = \langle T', M_1 \rangle, k_2 = \langle T', M_2 \rangle, k_3 = \langle T', B_2 \rangle$$

as parallel transport curvatures.

The Euler angles are introduced by Leonhard Euler to describe the orientation of a rigid body. Also, the Euler angles are used in robotics for speaking about the degrees of freedom of a wrist, electronic stability control in a smilar way, gun fire control systems require corrections to gun order angles to compensate for dect tilt, crystallographic texture can be described using Euler angles and texture analysis, the Euler angles provide the necessary mathematical depiction of the orientation of individual crystallites within a polycrystalline material, allowing for the quantitative description of the macroscopic material. [see more details in [18]], the most application is to describe aircraft attitudes see in Figure 2.

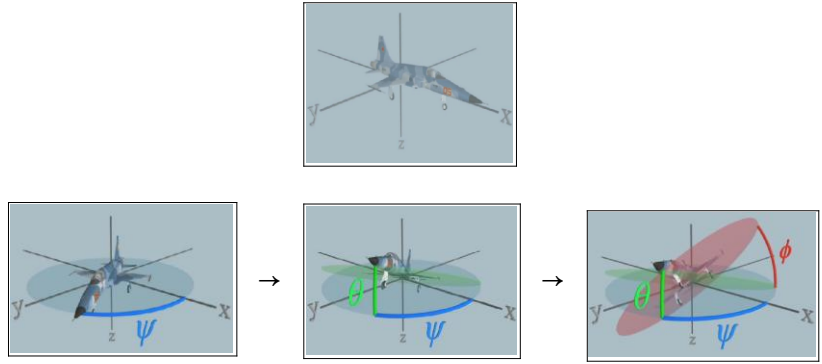

Figure 2: Euler angles of the aircraft movements.

Using Euler angles an arbitrary rotation matrix is given by

$$\begin{bmatrix} \cos\theta\cos\psi & -\cos\phi\sin\psi+\sin\phi\sin\theta\sin\psi & \sin\phi\sin\psi+\cos\phi\sin\theta\cos\psi \\ \cos\theta\sin\psi & \cos\phi\cos\psi+\sin\phi\sin\theta\sin\psi & -\sin\phi\cos\psi+\cos\phi\sin\theta\sin\psi \\ -\sin\theta & \sin\phi\cos\theta & \cos\phi\cos\theta \end{bmatrix}$$

where $\theta, \phi, \psi$ are Euler angles [18].

## 3. Parallel Transport Frame in $4-$dimensional Euclidean Space $\mathbb{E}^4$

In this section, we give parallel transport frame of a curve and we introduce the relations between the frame and Frenet frame of the curve in $4-$dimensional Euclidean space using the Euler angles. The relation which is well known in Euclidean $3-$space is generalized for the first time in $4-$dimensional Euclidean space.

**Theorem 1.** *Let $\{T, M_1, M_2, M_3\}$ be a Frenet frame along a unit speed curve $\alpha : I \subset \mathbb{R} \to \mathbb{E}^4$ and $\{T, M_1, M_2, M_3\}$ denotes the parallel transport frame of the curve $\alpha$. The relation may be expressed as*

$$T = T(s)$$
$$N = \cos\theta(s)\cos\psi(s)M_1 + (-\cos\phi(s)\sin\psi(s) + \sin\phi(s)\sin\theta(s)\cos\psi(s))M_2$$
$$\quad + (\sin\phi(s)\sin\psi(s) + \cos\phi(s)\sin\theta(s)\cos\psi(s))M_3$$
$$B_1 = \cos\theta(s)\sin\psi(s)M_1 + (\cos\phi(s)\cos\psi(s) + \sin\phi(s)\sin\theta(s)\sin\psi(s))M_2$$
$$\quad + (-\sin\phi(s)\cos\psi(s) + \cos\phi(s)\sin\theta(s)\sin\psi(s))M_3$$
$$B_2 = -\sin\theta(s)M_1 + \sin\phi(s)\cos\theta(s)M_2 + \cos\phi(s)\cos\theta(s)M_3$$

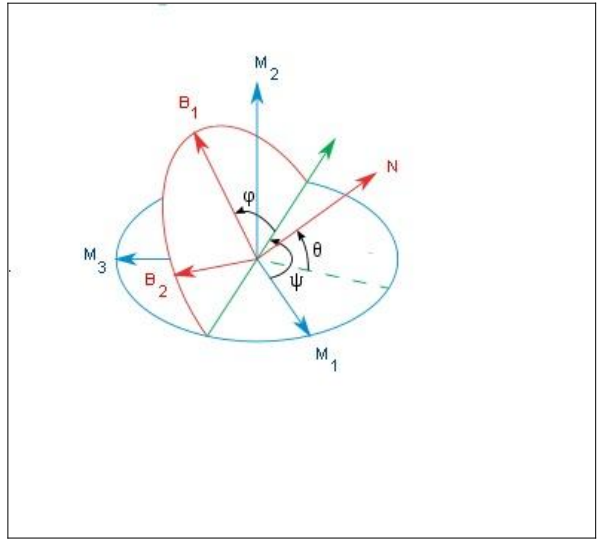

Figure 3: The relation of the Frenet and the parallel transport frame
by means of the rotation matrix.

*The alternative parallel frame equations are*

$$\begin{bmatrix} T' \\ M_1' \\ M_2' \\ M_3' \end{bmatrix} = \begin{bmatrix} 0 & k_1 & k_2 & k_3 \\ -k_1 & 0 & 0 & 0 \\ -k_2 & 0 & 0 & 0 \\ -k_3 & 0 & 0 & 0 \end{bmatrix} \begin{bmatrix} T \\ M_1 \\ M_2 \\ M_3 \end{bmatrix} \quad (3.1)$$

*where $k_1, k_2, k_3$ are principal curvature functions according to parallel transport frame of the curve $\alpha$ and their expression as follows*

$$k_1 = \kappa \cos\theta \cos\psi,$$
$$k_2 = \kappa(-\cos\phi \sin\psi + \sin\phi \sin\theta \cos\psi),$$
$$k_3 = \kappa(\sin\phi \sin\psi + \cos\phi \sin\theta \cos\psi)$$

*where* $\theta' = \dfrac{\sigma}{\sqrt{\kappa^2 + \tau^2}}$, $\psi' = -\tau - \sigma \dfrac{\sqrt{\sigma^2 - \theta'^2}}{\sqrt{\kappa^2 + \tau^2}}$, $\phi' = -\dfrac{\sqrt{\sigma^2 - \theta'^2}}{\cos\theta}$ *and the following equalities*

$$\kappa(s) = \sqrt{k_1^2 + k_2^2 + k_3^2},$$
$$\tau(s) = -\psi' + \phi' \sin\theta,$$
$$\sigma(s) = \frac{\theta'}{\sin\psi},$$
$$\phi' \cos\theta + \theta' \cot\psi = 0.$$

*hold.*

**Proof.** Let $\{T, M_1, M_2, M_3\}$ be a Frenet frame along a unit speed curve $\alpha : I \subset \mathbb{R} \to \mathbb{E}^4$ and $\{T, M_1, M_2, M_3\}$ denotes the parallel transport frame of the curve $\alpha$. The relation between Frenet frame and parallel transport frame as follows

$$T = T(s)$$
$$N = \cos\theta(s)\cos\psi(s)M_1 + (-\cos\phi(s)\sin\psi(s) + \sin\phi(s)\sin\theta(s)\cos\psi(s))M_2$$
$$+ (\sin\phi(s)\sin\psi(s) + \cos\phi(s)\sin\theta(s)\cos\psi(s))M_3$$
$$B_1 = \cos\theta(s)\sin\psi(s)M_1 + (\cos\phi(s)\cos\psi(s) + \sin\phi(s)\sin\theta(s)\sin\psi(s))M_2$$
$$+ (-\sin\phi(s)\cos\psi(s) + \cos\phi(s)\sin\theta(s)\sin\psi(s))M_3$$
$$B_2 = -\sin\theta(s)M_1 + \sin\phi(s)\cos\theta(s)M_2 + \cos\phi(s)\cos\theta(s)M_3$$

Differentiating the $M_1$, $M_2$, $M_3$ with respect to $s$ we get

$$M_1' = (-\kappa\cos\theta\cos\psi)T$$
$$+ (\theta'\sin\theta\cos\psi - \psi'\cos\theta\sin\psi - \tau\cos\theta\sin\psi)N$$
$$+ (-\theta'\sin\theta\sin\psi + \psi'\cos\theta\cos\psi + \sigma\sin\theta)B_1$$
$$+ (-\theta'\cos\theta + \sigma\cos\theta\sin\psi)B_2,$$

$$M_2' = -\kappa(-\cos\phi\sin\psi + \sin\phi\sin\theta\sin\psi)T$$
$$+ (\phi'\sin\phi\sin\psi - \psi'\cos\phi\cos\psi + \phi'\cos\phi\sin\theta\cos\psi$$
$$+ \theta'\sin\phi\cos\theta\cos\psi - \psi'\sin\phi\sin\theta\sin\psi - k_2(\cos\phi\cos\psi + \sin\phi\sin\theta\sin\psi))N$$
$$+ (-\phi'\sin\phi\cos\psi - \psi'\cos\phi\sin\psi + \phi'\cos\phi\sin\theta\sin\psi$$
$$+ \theta'\sin\phi\cos\theta\sin\psi + \psi'\sin\phi\sin\theta\cos\psi + \tau(-\cos\phi\sin\psi + \sin\phi\sin\theta\cos\psi)$$
$$- \sigma\sin\phi\cos\theta)B_1$$
$$+ (\phi'\cos\phi\cos\theta - \theta'\sin\phi\sin\theta$$
$$+ \sigma(\cos\phi\cos\psi + \sin\phi\sin\theta\sin\psi))B_2$$

$$M_3'(s) = -\kappa(\sin\phi\sin\psi + \cos\phi\sin\theta\cos\psi)T$$
$$+ (\phi'\cos\phi\sin\psi + \psi'\sin\phi\cos\psi - \phi'\sin\phi\sin\theta\cos\psi$$
$$+ \theta'\cos\phi\cos\theta\cos\psi - \psi'\cos\phi\sin\theta\sin\psi - \tau(-\sin\phi\cos\psi + \cos\phi\sin\theta\sin\psi))N$$
$$+ (-\phi'\cos\phi\cos\psi + \psi'\sin\phi\sin\psi - \phi'\sin\phi\sin\theta\sin\psi + \theta'\cos\phi\cos\theta\sin\psi$$
$$+ \psi'\cos\phi\sin\theta\cos\psi + \tau(\sin\phi\sin\psi + \cos\phi\sin\theta\cos\psi)$$
$$+ \sigma(\cos\phi\cos\theta))B_1$$
$$+ (-\phi'\sin\phi\cos\theta - \theta'\cos\phi\sin\theta$$
$$+ \sigma(-\sin\phi\cos\psi + \cos\phi\sin\theta\sin\psi))B_2.$$

Since the $M_1$, $M_2$ and $M_3$ are relatively parallel vector field, normal component of the $M_1'$,

$M_2'$ and $M_3'$ must be zero and the following equalities $\langle M_1', M_2 \rangle = 0$, $\langle M_1', M_3 \rangle = 0$ satisfy.
Also, if we consider that $k_1 = \langle T', M_1 \rangle$, $k_2 = \langle T', M_2 \rangle$, $k_3 = \langle T', M_3 \rangle$, we can easily see that

$$\kappa(s) = \sqrt{k_1^2 + k_2^2 + k_3^2},$$
$$\tau(s) = -\psi' + \phi' \sin \theta,$$
$$\sigma(s) = \frac{\theta'}{\sin \psi},$$
$$\phi' \cos \theta + \theta' \cot \psi = 0.$$

and

$$k_1 = \kappa \cos \theta \cos \psi,$$
$$k_2 = \kappa(-\cos \phi \sin \psi + \sin \phi \sin \theta \cos \psi),$$
$$k_3 = \kappa(\sin \phi \sin \psi + \cos \phi \sin \theta \cos \psi).$$

If we choose $\theta' = \frac{\sigma}{\sqrt{\kappa^2 + \tau^2}}$, then we obtain $\psi' = -\tau - \sigma \frac{\sqrt{\sigma^2 - \theta'^2}}{\sqrt{\kappa^2 + \tau^2}}$, $\phi' = -\frac{\sqrt{\sigma^2 - \theta'^2}}{\cos \theta}$. So, we have

$$\begin{bmatrix} T' \\ M_1' \\ M_2' \\ M_3' \end{bmatrix} = \begin{bmatrix} 0 & k_1 & k_2 & k_3 \\ -k_1 & 0 & 0 & 0 \\ -k_2 & 0 & 0 & 0 \\ -k_3 & 0 & 0 & 0 \end{bmatrix} \begin{bmatrix} T \\ M_1 \\ M_2 \\ M_3 \end{bmatrix},$$

which complete the proof.

**Corollary 1.** *If we consider that $\theta = \phi = 0$ then we get the Bishop frame in $\mathbb{E}^3$.*
Now, we give an example for the curve which has not a Frenet frame at some points but, it has parallel transport frame on these points.

**Example 1.** *Let $\alpha(s) = (\sin s, 2s + 1, 2s - 1, s)$ be a curve in Euclidean $4-$space . Since $\alpha''(0) = (0, 0, 0, 0)$, we can not calculate the Frenet frame vectors at the point $s = 0$. However, we can calculate using the parallel transport frame vectors as follows*
$$M_1 = \sin \psi B_1,$$
$$M_2 = \cos \phi \cos \psi B_1,$$
$$M_3 = -\sin \phi \cos \psi B_1$$
*where $\psi$ and $\phi$ constant angles.*

**Theorem 2.** *Let $\alpha : I \subset \mathbb{R} \to \mathbb{E}^4$ be a curve with nonzero curvatures $k_i$ $(i = 1, 2, 3)$ according to parallel transport frame in Euclidean $4-$space $\mathbb{E}^4$. Then $\alpha$ lies on a sphere if and only if $ak_1 + bk_2 + ck_3 + 1 = 0$ where $a$, $b$, and $c$ are non-zero constants.*

**Proof.** Let $\alpha$ lies on a sphere with center $P$ and radius $r$, then $\langle \alpha - P, \alpha - P \rangle = r^2$. Differentiating this equation with respect to $s$, it gives us $\langle T, \alpha - P \rangle = 0$, so $\alpha - P = aM_1 + bM_2 + cM_3$ for some function $a$, $b$, and $c$. Where

$a' = \langle \alpha - P, M_1 \rangle' = \langle T, M_1 \rangle + \langle -k_1 T, \alpha - P \rangle = 0$ so $a$ is a constant. Similarly, we can easily say that $b$ and $c$ are constant. Then differentiating the equation $\langle T, \alpha - P \rangle$ with respect to $s$, we get $\langle k_1 M_1 + k_2 M_2 + k_3 M_3, \alpha - P \rangle + \langle T, T \rangle = ak_1 + bk_2 + ck_3 + 1 = 0$. That is, between $k_1, k_2$ and $k_3$ have the linear relation such as
$$ak_1 + bk_2 + ck_3 + 1 = 0.$$
Moreover, $r^2 = \langle \alpha - P, \alpha - P \rangle = a^2 + b^2 + c^2 = \frac{1}{d^2}$ where $d$ is the distance of the plane $ax + by + cz + 1 = 0$ from the origin.
Conversely, suppose that the equation $ak_1 + bk_2 + ck_3 + 1 = 0$ holds. If $P$ is denoted by $P = \alpha - aM_1 - bM_2 - cM_3$ then differentiating the last equation we have $P' = T + (ak_1 + bk_2 + ck_3)T = 0$ so $P$ is constant. Similarly shows that $r^2 = \langle \alpha - P, \alpha - P \rangle$ is constant. So, $\alpha$ lies on a sphere with center $P$ and radius $r$.

**Example 2.** *Let* $\alpha(s) = (\sin \frac{s}{\sqrt{2}}, \cos \frac{s}{\sqrt{2}}, \frac{1}{\sqrt{2}} \sin s, \frac{1}{\sqrt{2}} \cos s)$ *be a curve in Euclidean* $4-space$ .
*According the Frenet frame there are lots of formulas for showing that this curve is a spherical curve. But the formulas have some disadvantages which were define the above chapters. Then we calculate curvature functions of the curve* $\alpha$ *according to parallel transport frame*
$$k_1 = 0, \ k_2 = -\cos\phi, \ k_3 = \sin\phi$$
*where* $\phi$ *is constant. The curve* $\alpha$ *satisfy the fallowing equation*
$$ak_1 + bk_2 + ck_3 + 1 = 0.$$
*Consequently, the curve* $\alpha$ *is a spherical curve. But, using the Frenet curvatures we can not show that* $\alpha$ *is a spherical curve. Because* $\alpha$ *has a zero torsion.*

## 4. Normal, rectifying and osculating curves according to parallel transport frame

In this chapter, we define normal, rectifying and osculating curves according to parallel transport frame and obtain some characterizations for such curves.

### 4.1. Normal curves according to parallel transport frame

The normal space according to parallel transport frame of the curve $\alpha$ as the orthogonal complement $T^\perp$ of its tangent vector field $T$. Hence the normal space is given by $T^\perp = \{X \in \mathbb{E}^4 | \langle X, T \rangle = 0\}$. Normal curves are defined in [12] as a curve whose position vector always lies in its normal space. Consequently, the position vector of the normal curve $\alpha$ with parallel transport vector fields $M_1, M_2$ and $M_3$ satisfies the equation
$$\alpha(s) = \lambda(s) M_1(s) + \mu(s) M_2(s) + \upsilon(s) M_3(s), \qquad (4.1)$$
where $\lambda(s), \mu(s)$ and $\upsilon(s)$ differentiable functions.
Then we have the following theorem which characterize normal curves according to parallel transport frame.

**Theorem 3.** *Let $\alpha(s)$ be a unit speed curve in Euclidean $4-$space with parallel transport vector fields $T, M_1, M_2, M_3$ and its curvature functions $k_1, k_2, k_3$ of the curve $\alpha$. Then the curve $\alpha$ is a normal curve if and only if $\alpha$ is a spherical curve.*

**Proof.** Let us first assume that $\alpha$ is a normal curve. Then its position vector satisfies the Equation (4.1). By taking the derivative of Eq. (4.1) with respect to $s$ and using the parallel transport frame in Eq. (3.1), we obtain

$$T(s) = -(\lambda(s)k_1(s) + \mu(s)k_2(s) + \upsilon(s)k_3(s))T(s)$$
$$+ \lambda'(s)M_1(s) + \mu'(s)M_2(s) + \upsilon'(s)M_3(s)$$

and so

$$\lambda'(s) = 0$$
$$\mu'(s) = 0,$$
$$\upsilon'(s) = 0.$$

From last equations we find

$$\lambda(s)k_1(s) + \mu(s)k_2(s) + \upsilon(s)k_3(s) + 1 = 0,$$

where $\lambda(s), \mu(s)$ and $\upsilon(s)$ are non zero constant functions. From Theorem 2 the curve $\alpha$ is a spherical curve.

Conversely, suppose that the curve $\alpha$ is a spherical curve. Let us consider the vector $m \in \mathbb{E}^4$ given by

$$m(s) = \alpha(s) - (\lambda(s)M_1(s) + \mu(s)M_2(s) + \upsilon(s)M_3(s)) \qquad (4.2)$$

Differentiating (4.2) with respect to $s$ and applying (3.1), we get

$$m'(s) = (1 + \lambda(s)k_1(s) + \mu(s)k_2(s) + \upsilon(s)k_3(s))T(s)$$
$$+ \lambda'(s)M_1(s) + \mu'(s)M_2(s) + \upsilon'(s)M_3(s)$$

Since the curve $\alpha$ is a spherical curve and $\lambda(s), \mu(s)$ and $\upsilon(s)$ are non zero constant we have $m'(s) = 0$. Consequently, $\alpha$ is a normal curve.

### 4.2. Rectifying curves according to parallel transport frame

The rectifying space according to parallel transport frame of the curve $\alpha$ as the orthogonal complement $M_1^\perp$ of its first normal vector field $M_1$. Hence the normal space is given by $M_1^\perp = \{X \in \mathbb{E}^4 | \langle X, M_1 \rangle = 0\}$. Rectifying curves are defined in [il]. as a curve whose position vector always lies in its rectifying space. Consequently, the position vector of the rectifying curve $\alpha$ with respect to parallel transport vector fields $T$, $M_2$ and $M_3$ satisfies the equation

$$\alpha(s) = c_1(s)T(s) + c_2(s)M_2(s) + c_3(s)M_3(s), \qquad (4.3)$$

where $c_1(s)$, $c_2(s)$ and $c_3(s)$ differentiable functions.
Then we have the following theorem which characterize rectifying curves according to parallel transport frame.

**Theorem 4.** *Let* $\alpha : I \subset \mathrm{R} \to \mathbb{E}^4$ *be a curve with nonzero curvatures* $k_i$ $(i=1,2,3)$ *according to parallel transport frame in* $\mathbb{E}^4$. *Then* $\alpha$ *is a rectifying curve if and only if* $c_2 k_2 + c_3 k_3 + 1 = 0$, $c_2, c_3 \in \mathrm{R}$.

**Proof.** Let us first assume that $\alpha$ be a rectifying curve. With the help of the definition of rectifying curve, the position vectors of the curve $\alpha$ satisfies
$$\alpha(s) = c_1 T(s) + c_2 M_2(s) + c_3 M_3(s). \tag{4.4}$$
Differentiating the last equation with respect to $s$ and using the parallel transport frame formulas we get
$$T = (c_1' - c_2 k_2 - c_3 k_3)T + c_1 k_1 M_1 + (c_2' + c_1 k_2)M_2 + (c_3' + c_1 k_3)M_3 \tag{4.5}$$
and therefore
$$c_1 = 0$$
$$c_2 = \text{const.}$$
$$c_3 = \text{const.} \tag{4.6}$$
Then using the last equation, we can easily find that the curvatures $k_1, k_2$ and $k_3$ satisfies the equation
$$c_2 k_2 + c_3 k_3 + 1 = 0 \tag{4.7}$$
where $c_2$ and $c_3$ are constants.
Conversely, suppose that the curvatures $k_1, k_2$ and $k_3$ satisfies the equation $c_2 k_2 + c_3 k_3 + 1 = 0$, $c_2, c_3 \in \mathrm{R}$. Let us consider the vector $X \in \mathbb{E}^4$ given by
$$X(s) = \alpha(s) - c_2 M_2(s) + c_3 M_3(s) \tag{4.8}$$
then using the Eq. (4.6) and (4.7) we can easily see that $X'(s) = 0$, that is, $X$ is constant vector. So, $\alpha$ is a rectifying curve. This completes the proof.

**Corollary 2.** *Let* $\alpha : I \subset \mathrm{R} \to \mathbb{E}^4$ *be a curve with nonzero curvatures* $k_i$ $(i=1,2,3)$ *according to parallel transport frame in* $\mathbb{E}^4$. *Then* $\alpha$ *is a rectifying curve if and only if*
$$(\frac{k_3}{k_2})' \frac{k_2^2}{k_3'} = \text{const. or } (\frac{k_3}{k_1})' \frac{k_1^2}{k_3'} = \text{const.}$$

**Proof.** The proof is clear from Eq. (4.7).

### 4.3 Osculating curves according to parallel transport frame

The osculating space according to parallel transport frame of the curve $\alpha$ as the orthogonal complement $M_2^\perp$ of its second binormal vector of the parallel transport frame. Hence the osculating space is given by $M_2^{\perp} = \{X \in \mathbb{E}^4 | \langle X, M_2 \rangle = 0\}$. Osculating curves are defined in [11]. as a curve whose position vector always lies in its osculating space. Consequently, the position vector of the osculating curve $\alpha$ with respect to parallel transport vector fields $T$, $M_1$ and $M_3$ satisfies the equation

$$\alpha(s) = \lambda_1(s)T(s) + \lambda_2(s)M_1(s) + \lambda_3(s)M_3(s), \tag{4.9}$$

where $\lambda_1(s)$, $\lambda_2(s)$ and $\lambda_3(s)$ differentiable functions.

Then we have the following theorem which characterize osculating curves according to parallel transport frame.

**Theorem 5.** *Let* $\alpha : I \subset \mathbb{R} \to \mathbb{E}^4$ *be a curve with nonzero curvatures* $k_i$ $(i = 1, 2, 3)$ *according to parallel transport frame in* $\mathbb{E}^4$. *Then* $\alpha$ *is a osculating curve if and only if* $\lambda_2 k_1 + \lambda_3 k_3 + 1 = 0$, $\lambda_2, \lambda_3 \in \mathbb{R}$.

**Proof.** Let us first assume that $\alpha$ be a osculating curve. With the help of the definition of osculating curve, the position vectors of the curve $\alpha$ satisfies

$$\alpha(s) = \lambda_1 T(s) + \lambda_2 M_1(s) + \lambda_3 M_3(s) \tag{4.10}$$

Differentiating the last equation with respect to $s$ and using the parallel transport equations we get

$$T = (\lambda_1' - \lambda_2 k_1 - \lambda_3 k_3)T + (\lambda_2' + \lambda_1 k_2)M_2 + \lambda_1 k_1 M_1 + (\lambda_3' + \lambda_1 k_3)M_3 \tag{4.11}$$

and therefore

$$\lambda_1 = 0$$
$$\lambda_2 = const.$$
$$\lambda_3 = const. \tag{4.12}$$

Then using the last equation, we can easily find that the curvatures $k_1, k_2$ and $k_3$ satisfies the equation

$$\lambda_2 k_1 + \lambda_3 k_3 + 1 = 0 \tag{4.13}$$

where $\lambda_2$ and $\lambda_3$ are constants.

Conversely, suppose that the curvatures $k_1, k_2$ and $k_3$ satisfies the equation $\lambda_2 k_1 + \lambda_3 k_3 + 1 = 0$, $\lambda_2, \lambda_3 \in \mathbb{R}$. Let us consider the vector $X \in \mathbb{E}^4$ given by

$$X(s) = \alpha(s) - \lambda_2 M_1(s) + \lambda_3 M_3(s) \tag{4.14}$$

then using the Eq. (4.12) and (4.13) we can easily see that $X'(s) = 0$, that is, $X$ is constant vector. So, $\alpha$ is a osculating curve. This completes the proof.

**Corollary 3.** *Let* $\alpha : I \subset \mathbb{R} \to \mathbb{E}^4$ *be a curve with nonzero curvatures* $k_i$ $(i = 1, 2, 3)$ *according to parallel transport frame in* $\mathbb{E}^4$. *Then* $\alpha$ *is a rectifying curve if and only if*

$$(\frac{k_3}{k_1})' \frac{k_1^2}{k_3'} = const. \text{ or } (\frac{k_3}{k_1})' \frac{k_1^2}{k_2'} = const.$$

**Proof.** The proof is clear from Eq. (4.13).

Fatma GÖKÇELİK
Department of Mathematics, Faculty of Science, University of Ankara Tandogan, Ankara, TURKEY
*E-mail :* fgokcelik@ankara.edu.tr
Zehra BOZKURT
Department of Mathematics, Faculty of Science, University of Ankara Tandogan, Ankara, TURKEY
*E-mail :* zbozkurt@ankara.edu.tr
İsmail GÖK
Department of Mathematics, Faculty of Science, University of Ankara Tandogan, Ankara, TURKEY
*E-mail :* igok@science.ankara.edu.tr
F. Nejat EKMEKCİ
Department of Mathematics, Faculty of Science, University of Ankara Tandogan, Ankara, TURKEY
*E-mail :* ekmekci@science.ankara.edu.tr
Yusuf YAYLI
Department of Mathematics, Faculty of Science, University of Ankara Tandogan, Ankara, TURKEY
*E-mail :* yyayli@science.ankara.edu.tr